\documentclass{amsart}
\usepackage[utf8]{inputenc}

\usepackage{graphics}

\usepackage[T1]{fontenc}    % use 8-bit T1 fonts

\usepackage{amsthm}
\usepackage{amsbsy,amsmath,amssymb,amscd,amsfonts}
\usepackage{hyperref}
\usepackage{url}            % simple URL typesetting
\usepackage{booktabs}       % professional-quality tables
\usepackage{nicefrac}       % compact symbols for 1/2, etc.
\usepackage{microtype}      % microtypography

\usepackage{graphicx,float,latexsym,color}
\usepackage[font={scriptsize,it}]{caption}
\usepackage{subcaption}

\usepackage{makecell}

\usepackage[dvipsnames]{xcolor}

\hypersetup{
    pdftoolbar=true,        % show Acrobat’s toolbar?
    pdfmenubar=true,        % show Acrobat’s menu?
    pdffitwindow=false,     % window fit to page when opened
    pdfstartview={FitH},    % fits the width of 
    colorlinks=true,       % false: boxed links; true: colored links
    linkcolor=OliveGreen,          % color of internal links (change box color with linkbordercolor)
    citecolor=blue,        % color of links to bibliography
    filecolor=black,      % color of file links
    urlcolor=red           % color of external links
}

\usepackage{lineno}

\renewcommand{\P}{\mathcal{P}}

\title[Eighty New Billiard Invariants]{Eighty New Invariants\\in the Elliptic Billiard}
\author{Dan Reznik}
\address{Data Science Consulting\\Rio de Janeiro, Brazil}
\email{dreznik@gmail.com}
\author{Ronaldo Garcia}
\address{Math \& Statistics Institute\\Federal University of Goiás\\Goiânia, Brazil}
\email{ragarcia@ufg.br}
\author{Jair Koiller}
\address{Federal University of Rio de Janeiro\\Rio de Janeiro, Brazil}
\email{jairkoiller@gmail.com}
\date{October, 2020}

\newcommand{\Mod}[1]{\ (\mathrm{mod}\ #1)}

\begin{document}

\maketitle

\begin{abstract}
    We introduce several-dozen experimentally-found invariants of Poncelet N-periodics in the confocal ellipse pair (Elliptic Billiard). Recall this family is fully defined by two integrals of motion (linear and angular momentum), so any ``new'' invariants are dependent upon them. Nevertheless, proving them may require sophisticated methods. We reference some two-dozen proofs already contributed. We hope this article will motivate contributions for those still lacking proof.
    

\end{abstract}

\section{Introduction}
\begin{figure}[H]
    \centering
    \includegraphics[width=.5\textwidth]{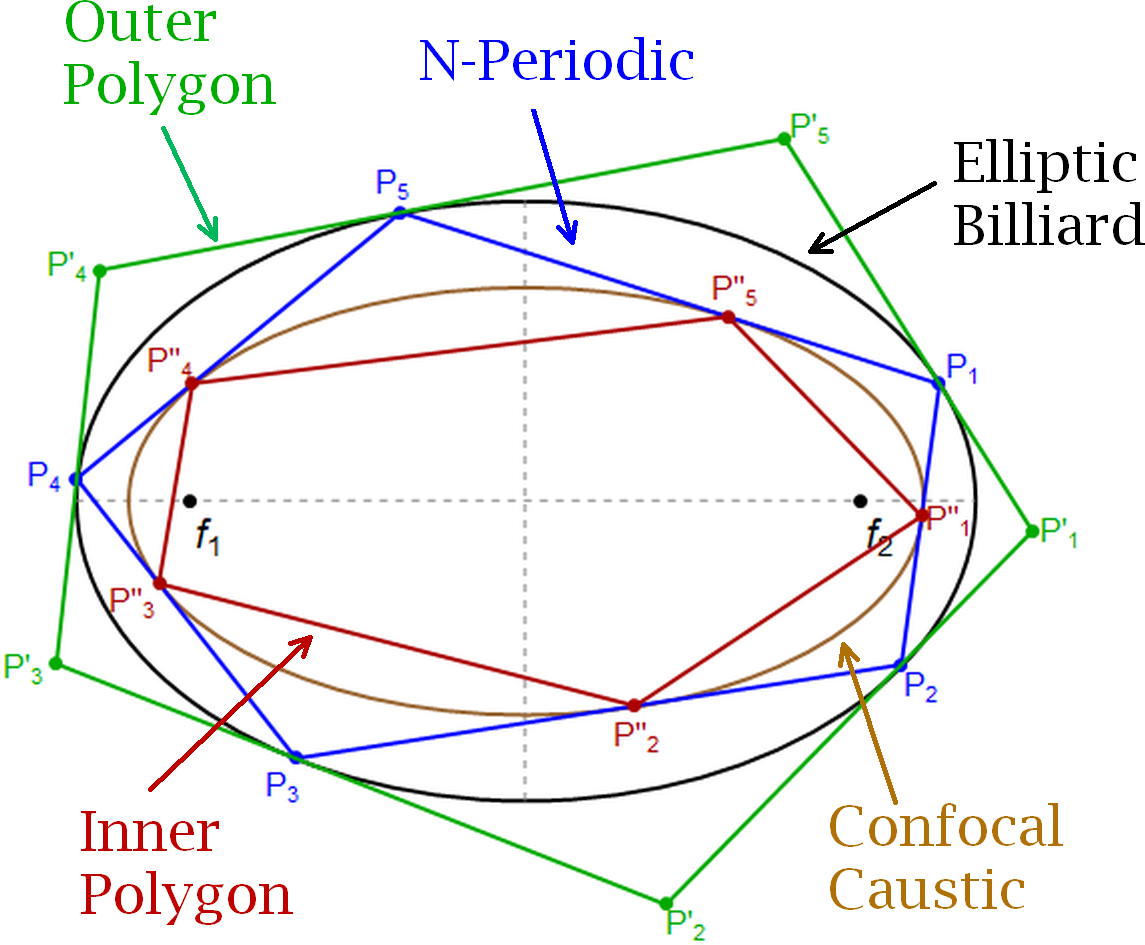}
    \caption{A (i) 5-periodic (vertices $P_i$) is shown inscribed in a confocal ellipse pair (billiard and caustic). Also shown is (ii) the outer polygon with vertices $P_i'$ tangent to the outer ellipse at the $N$-periodic vertices, and (iii) the inner polygon whose vertices $P_i''$ are at the points of contact of the $N$-periodic with the caustic.}
    \label{fig:eb}
\end{figure}

 The Elliptic Billiard (EB) is a special case of Poncelet's Porism \cite{dragovic2014-poncelet}, where the conic pair are two confocal ellipses; it therefore admits a 1d family of $N$-periodic trajectories \cite{sergei91,izmestiev2017-ivory,dragovic11} which at every vertex are bisected by the normals of the outer ellipse in the pair (hence the term ``billiard''); see Figure~\ref{fig:eb}.
 
The EB is an {\em integrable} system (in fact it is conjectured as the only integrable planar billiard \cite{kaloshin2018}).
 Integrability implies invariant perimeter $L$; a second classic invariant is Joachimsthal's constant $J$, which is simply a statement that all trajectory segments are tangent to the confocal caustic \cite{sergei91,rozikov2018}.

Continuing our work on properties of N-periodics in the EB \cite{reznik2019-intelligencer,garcia2020-new-properties}, here we introduce dozens of new invariants  detected experimentally. These involve distances, areas, angles and centers of mass of $N$-periodics and several derived polygons defined below. Some invariants depend on the parity of $N$, while others on positional constraints.

Note that since the N-periodics in the EB are fully defined by $L,J$, any ``new'' invariants listed here or elsewhere must be ultimately dependent upon said quantities. Nevertheless, proving a specific functional dependence may require sophisticated techniques. Several proofs have already been contributed and are referenced below. We hope to motivate more contributions and/or new discoveries.

Admittedly, the number of possible invariants is infinite as one may select any functional combination of $L,J$. Our selection criterion can be loosely defined as {\em any quantity constant over the N-periodic family and/or derived objects, which is an elementary function of lengths, angles, areas, etc}.

This article is organized as follows: preliminary definitions are given in Section~\ref{sec:prelim}. Invariants are introduced in Section~\ref{sec:invariants}, in several clusters, involving: (i) lengths, areas, and angles of $N$-periodics and associated polygons; (ii) pedal polygons to N-periodics and (iii) their outer polygons; (iv)  antipedal polygons (defined below); (v) area-ratios related to the Steiner  curvature centroid \cite{steiner1838}; (vi) pairs of pedal polygons; (vii) area-ratios of evolute polygons \cite{fuchs2017}; (viii) focus-inversive objects and (ix) pairs of focus-inversive objects.

Details about our experimental toolbox are covered in Section~\ref{sec:method}. All symbols used in this article are listed on Table~\ref{tab:symbols} in Appendix~\ref{app:symbols}.

We encourage the reader to watch the videos included in Section~\ref{sec:videos} which provide more insight into invariant phenomena. 

\subsection*{Related Work}
In a companion article \cite{garcia2020-invariants} we derive explicit expressions for some of the invariants listed herein for certain ``low'' $N$, e.g., 3--6. Methods for obtaining N-periodic trajectories based on Cayley's condition are surveyed in \cite{dragovic2006-survey,chang1988}. A few explicit expressions for the caustic parameter (for $N=3,4,6,8$) appear in \cite{ros2014-ellipsoids}.

\section{Preliminaries}
\label{sec:prelim}
Let the EB have center $O$, semi-axes $a>b>0$, and foci $f_1,f_2$ at $[{\pm}\sqrt{a^2-b^2},0]$. Let $a'',b''$ denote the major, minor semi-axes of the confocal caustic, whose values are given by a method due to Cayley \cite{dragovic11}, though we obtain them numerically, see Section~\ref{sec:method}.

As mentioned above, the perimeter $L$ is invariant for a given $N$-periodic family, as is Joachmisthal's constant $J=\left<\mathcal{A}x,v\right>$, where $x$ is a bounce point (called $P_i$ above), $v$ is the unit velocity vector $(P_i-P_{i-1})/||.||$, $\left<.\right>$ stands for dot product, and \cite{sergei91}:

\[ \mathcal{A}=\mbox{diag}\left[1/a^2,1/b^2\right] \]

Hellmuth Stachel derived \cite{stachel2020-private-joachimsthal} an elegant expression for Joahmisthal's constant $J$ in terms of the axes of the EB and its caustic:

\begin{equation*}
    J = \frac{\sqrt{a^2 - a''^2}}{{a}{b}}
\end{equation*}

\noindent Note: holding $a$ constant, for each $N$, $a''$ and therefore $J$ assume a distinct value.

Let a polygon have vertices $W_i,i=1,...,N$. In this paper all polygon areas are {\em signed}, i.e., obtained from a sum of cross-products \cite{preparata1988}:

\begin{equation}
\mbox{S}=\frac{1}{2}\sum_{i=1}^N{W_i{\times}W_{i+1}}
\label{eqn:signed-a}
\end{equation}

Let $W_i=(x_i,y_i)$, then $W_i{\times}W_{i+1}=(x_i\,y_{i+1}-x_{i+1}\,y_i)$.\\

The area centroid $\overline{W}$ of a polygon is given by \cite{preparata1988}:

\begin{equation}
\overline{W}=\frac{1}{6S}\sum_{i=1}^N(W_i{\times}W_{i+1})(W_i+W_{i+1})
\label{eqn:area-centroid}
\end{equation}

The curvature $\kappa$ of the ellipse at point $(x,y)$ at distance $d_1,d_2$ to the foci is given by \cite[Ellipse]{mw}:

\begin{equation}
\kappa = \frac{1}{a^2 b^2} \left(\frac{x^2}{a^4}+\frac{y^2}{b^4}\right)^{-3/2} = a b (d_1 d_2)^{-3/2} = \left( \kappa_a d_1 d_2\right)^{-3/2}
\label{eqn:curv}
\end{equation}

\noindent Where $\kappa_a=(a b)^{-2/3}$ is the constant affine curvature of the ellipse \cite{nomizu1994-affine}.

Given a polygon with vertices $W_i$ and angles $\theta_i$, its Steiner Centroid of Curvature\footnote{J. Steiner (following a similar result by J. Sturm in 1823 for triangles) proved in 1825 that the area of pedal polygons of a polygon $W$ with respect to points on any given circumference centered on $K$ 
\cite{steiner1838} is invariant.} $K$ is given by \cite[p. 22]{steiner1838}:

\begin{equation}
K=
\frac{\sum_{i=1}^N{\rho_i}R_i}{\sum{\rho_i}},\,\,\,
\mbox{with}\,\rho_i=\sin(2\theta_i)
\label{eqn:steiner}
\end{equation}

\section{Invariants}
\label{sec:invariants}
In this section we present the invariants found so far in several tables. Each invariant is given an identifier $k_n$ where the first digit of $n$ refers to a cluster of invariants; see Table~\ref{tab:inv-group}.

\begin{table}[H]
%\small
\begin{tabular}{|c|l|l|}
\hline
range & invariant group & total \\
\hline
$k_{101}$--$k_{121}$ & Distances, area, angles, curvature & 21 \\
$k_{201}$--$k_{205}$ & N-Periodic Pedal polygons & 5 \\
$k_{301}$--$k_{307}$ & Outer pedal polygon & 7 \\
$k_{401}$--$k_{407}$ & Antipedal polygon & 7 \\
$k_{501}$--$k_{503}$ & Steiner curvature centroid & 3 \\
$k_{601}$--$k_{610}$ & Pairs of pedal polygons wrt. foci & 8 \\
$k_{701}$--$k_{703}$ & Evolute polygons & 3 \\
$k_{801}$--$k_{818}$ & Inversive objects & 18 \\
$k_{901}$--$k_{908}$ & Pairs of Inversive objects & 8 \\
\hline
\multicolumn{2}{r|}{\textbf{total:}} & \textbf{82} \\
\cline{3-3}
\end{tabular}
\caption{Numbering scheme for the invariants currently listed in this article.}
\label{tab:inv-group}
\end{table}

On the invariant tables below, column ``invariant'' provides an expression for the conserved quantity; column ``value'' provides a closed-form expression for the invariant (when available) in terms of the fundamental constants, or a `?' when not available (note that the invariant may already have been proved but no closed-form expression has yet been found); column ``which N'' specifies whether the invariant only holds for certain $N$ (even, odd, etc.); column ``date'' specifies the month and year (mm/yy) when the invariant was first experimentally detected. Column ``proven'' references available proofs if already communicated and/or published, else it displays a `?'.

\subsection{Basic Invariants}
\label{sec:basic}
Invariants involving angles and areas of N-periodics and its tangential and internal polygons are shown on Table~\ref{tab:invariants-general}. There $\theta_i$,$A$ (resp.~$\theta_i'$,$A'$) are angles, area of an N-periodic (resp.~outer polygon to the N-periodic). $A''$ is the area of the internal polygon (where orbit touches caustic), see Figure~\ref{fig:eb}. All sums/products go from $i=1$ to $N$. $k_{101},k_{102},k_{103}$ originally studied in \cite{reznik2019-intelligencer}. $l_i$ and $r_i$ denote $|P''_i-P_i|$ and $|P_{i+1}-P''_i|$,  respectively and $d_{j,i}=|P_i-f_j|$. $\kappa_i$ denotes the curvature of the EB at $P_i$ \eqref{eqn:curv}. $\alpha_{j,i}$ denotes the angle $P_i f_j P_{i+1}$.

%\begin{figure}
%    \centering
%    \includegraphics[width=\textwidth]{pics_pic_three_caustics.eps}
%    \caption{From left to right: $N$-periodic orbits (blue) for $N=4,5,6$, and their confocal caustic (brown). The {\em outer} polygon (green) is tangent to the EB (black) at every orbit vertex. Also shown is the {\em inner} polygon (red) whose vertices are the orbit's tangency points with the caustic.}
%    \label{fig:two-caustics}
%\end{figure}

\begin{table}
\begin{tabular}{|c|l|c|l|l|c|}
\hline
code & invariant & value & which N & date & proven \\
\hline
$k_{101}$ & $\sum{\cos\theta_i}$ & $JL-N$ & all & 4/19 & \cite{akopyan2020-invariants,bialy2020-invariants} \\
$k_{102}$ & $\prod{\cos\theta'_i}$ & ? & all & 5/19 & \cite{akopyan2020-invariants,bialy2020-invariants} \\
$k_{103}$ & $A'/A$ & ? & odd & 8/19 & \cite{akopyan2020-invariants,caliz2020-area-product} \\
$k_{104}$ & $\sum{\cos(2\theta'_i)}$ & ? & all & 1/20 & \cite{akopyan19_private_half_sines} \\ % \cite[Thm 6]{akopyan2020-invariants} \\
$k_{105}$ & $\prod{\sin(\theta_i/2)}$ & ? & odd & 1/20 &  \cite{akopyan19_private_half_sines} \\ % \cite[Thm 6]{akopyan2020-invariants} \\
$k_{106}$ & $A'A$ & ? & even & 1/20 & \cite{caliz2020-area-product} \\
$k_{107}$ & $k_{103} k_{105}$ & ? & ${\equiv}\,0 \Mod{4}$ & 1/20 & ? \\
$k_{108}$ & $k_{103} / k_{105}$ & ? & ${\equiv}\,2 \Mod{4}$ & 1/20 & ? \\
$k_{109}$ & $A/A''$ & $k_{103}$ & odd & 1/20 & ? \\
$k_{110}$ & $A\,A''$ & ? & even & 1/20 & ? \\
$k_{111}$ & $A'\,A''$ & ? & even & 1/20 & ? \\
$k_{112}$ & $A'\,A''/A^2$ & 1 & odd & 1/20 & \cite{akopyan2020-private-affine} \\
$k_{113}$ & $A'/A''$ & $[ab/({a''}{b''})]^2$ & all & 1/20 & \cite{stachel2020-private-k113-k116} \\
$k_{114}$ & $\prod{d_{1,i}}$ & ? & ${\equiv}\,2 \Mod{4}$ &  4/20 & ? \\
$k_{115}$ & $\prod|P_i'-f_1|$ & ? & ${\equiv}\,0 \Mod{4}$ & 4/20 & ? \\
${^\star}{k_{116}}$ & $\prod{l_i}/\prod{r_i}$ & 1 & all & 5/20 & \cite{stachel2020-private-k113-k116} \\
${^\star}{k_{117}}$ & $\prod{l_i},\;\prod{r_i}$ & ? & even & 5/20 & ?\\
${^\star}{k_{118}}$ & $\sum{l_i},\;\sum{r_i}$ & $L/2$ & odd & 8/20 & ?\\
$^{\dagger}{k_{119}}$ & $\sum{\kappa_i^{2/3}}$ & $L/[2 J (a b)^{4/3}]$ & all & 10/20 & \cite{stachel2020-private-curvature} \\
${^\ddag}k_{120}$ & $\sum\cos\alpha_{1,i}$& ? & all & 10/20 & ? \\
$k_{121}$ & $\sum{d_{1,i}}$ & ? & even & 10/20 & {\scriptsize symmetry} \\
\hline
\end{tabular}
\caption{Distance, area, and angle invariants displayed by the N-periodic, its outer and/or inner polygon. ${^\star}{k_i},i=116,117,118$ were discovered by Hellmuth Stachel. ${^\dagger}{k_{119}}$ was co-discovered with Pedro Roitman \cite{roitman2020-private} and is equivalent to $k_{902}$. ${^\ddag}{k_{120}}$ was suggested by A. Akopyan.}
\label{tab:invariants-general}
\end{table}

\subsection{Pedal Polygons}
\label{sec:pedal}
Tables~\ref{tab:inv-ped} and \ref{tab:inv-ped-tang} describe invariants found for the {\em pedal polygons} of N-periodics and the outer polygon, see Figure~\ref{fig:pic-pedals}.

\subsection{Pedals with respect to N-periodic}

\begin{figure}[H]
    \centering
    \includegraphics[width=\textwidth]{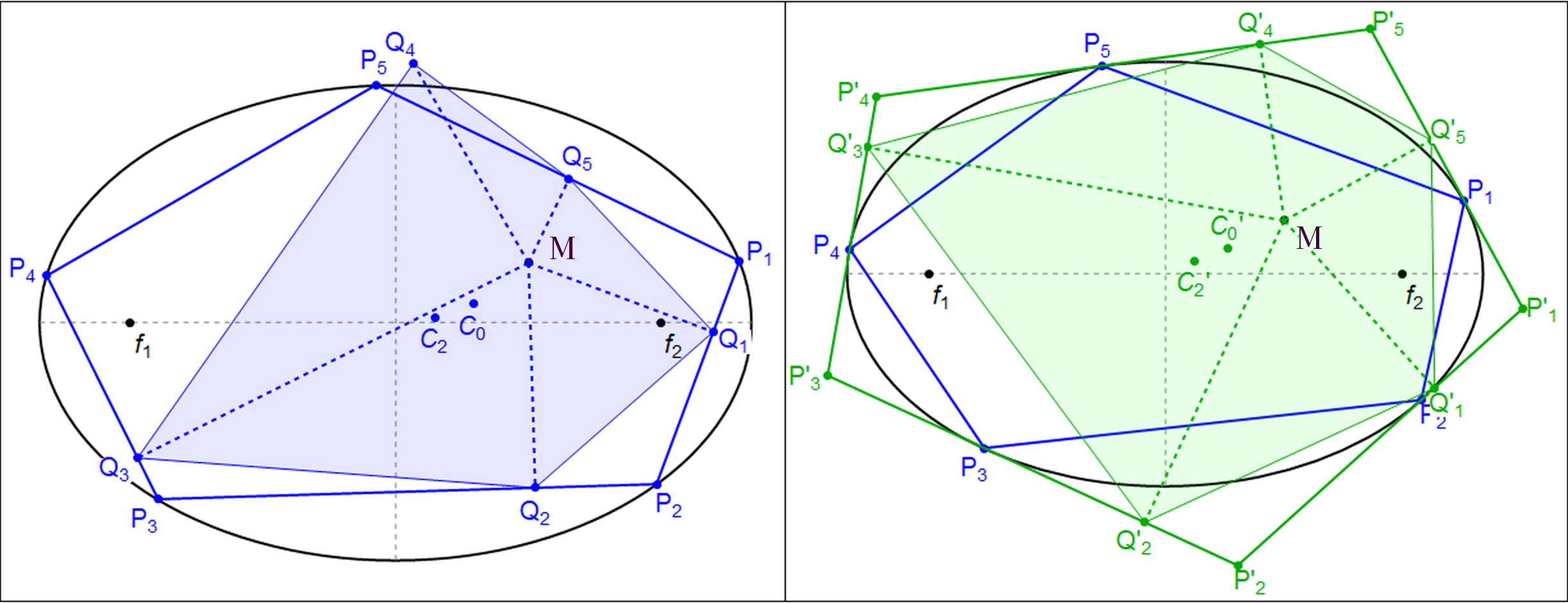}
    \caption{Left (resp. right): Pedal polygons for $N=5$ from a point $m$ with respect to the $N$-periodic (resp. its outer polygon). Vertex and area centroids $C_0,C_2$ are also shown. See Videos \cite[PL\#01,02,03]{dsr_playlist_2020_invariants}.}
    \label{fig:pic-pedals}
\end{figure}

Let $Q_i$ be the feet of perpendiculars dropped from a point $M$ onto the sides of the $N$-periodic. 
Let $A_m$ denote the area of the polygon formed by the $Q_i$, Figure~\ref{fig:pic-pedals}. Let $\phi_i$ denote the angle between two consecutive perpendiculars $Q_i-M$ and $Q_{i+1}-M$. Table~\ref{tab:inv-ped} lists invariants so far observed for these quantities.

\begin{table}[H]
\begin{tabular}{|c|l|c|l|c|l|c|}
\hline
code & invariant & value & which N & M & date & proven \\
\hline
${^\dagger}k_{201}$ & $|Q_i-O|$ & $a''$ & all & $f_1,f_2$ & 4/20 & \cite{akopyan20_private_cyclic_feet} \\
% subsumed by k_{20}
%$k_{15}$ & $\sum{d_i'^2}=\sum{e_i'^2}$ & ? & all & & 4/20 & ? \\
$k_{202,a}$ & $\prod|Q_i-M|$ & $(b'')^N$ & even & $f_1,f_2$ & 4/20 & \cite{bialy2020-invariants} \\
$k_{202,b}$ & $\prod|Q_i-M|$ & $(a''b'')^{N/2}$ & ${\equiv}\,0 \Mod{4}$ & O & 4/20 & \cite{bialy2020-invariants} \\
$k_{203,a}$ & $A\,A_m$ & ? & ${\equiv}\,0 \Mod{4}$ & all & 4/20 & ? \\
$k_{203,b}$ & $A\,A_m$ & ? & ${\not\equiv}\,2 \Mod{4}$ & O & 4/20 & ? \\
% subsumed by k_{20}
%$k_{18}$ & $A/A_m$ & ? & ${\equiv}\,2 \Mod{4}$ & O & 4/20 & ? \\
%$k_{17}$ & $A\,A_m$ & ? & all & O & 4/20 & ? \\
$k_{204}$ & $A/A_m$ & ? & ${\equiv}\,2 \Mod{4}$ & all & 4/20 & ? \\
$k_{205}$ & $\sum\cos\phi_i$ & ? & all & all & 4/20 & \cite{akopyan19_private_pedal_angles} \\
\hline
\end{tabular}
\caption{Invariants of pedal polygon with respect to N-Periodic sides. $^\dagger$ $k_{201}$ means the locus of the vertices of a pedal with respect to a focus is a circle.}
\label{tab:inv-ped}
\end{table}

\subsection{Pedals with respect to the Outer Polygon}

Let $Q_i'$ be the feet of perpendiculars dropped from a point $M$ onto the outer polygon. Let $\phi_i'$ denote the angle between two consecutive perpendiculars $Q'_i-M$ and $Q'_{i+1}-M$. Let $A_m'$ denote the area of the polygon formed by the $Q_i'$.

In the spirit of \cite{sergei2016-com} we also analyze centers of mass: $C_0'=\sum_i{Q_i'}/N$ is the vertex centroid, and the area centroid $C_2'$ of the polygon defined by the $Q_i'$ \eqref{eqn:area-centroid}. Table~\ref{tab:inv-ped-tang} lists invariants so far observed for these quantities.

\begin{table}[H]
\begin{tabular}{|c|l|c|l|c|l|c|}
\hline
code & invariant & value & which N & M & date & proven \\
\hline
${^\dagger}k_{301}$ & $|Q_i'-O|$ & $a$ & all & $f_1,f_2$ & 4/20 & \cite{akopyan20_private_cyclic_feet} \\
$k_{302}$ & 
$\sum{|Q_i'-M|^2}$ & ? & all & all & 4/20 & \cite{bialy2020-invariants} \\
$k_{303,a}$ & $A'\,A_m'$ & ? & ${\equiv}\,2 \Mod{4}$ &  all & 4/20 & ? \\
$k_{303,b}$ & $A'\,A_m'$ & ? & ${\not\equiv}\,0 \Mod{4}$ & O & 4/20 & ? \\
$k_{304}$ & $A'/A_m'$ & ? & ${\equiv}\,0 \Mod{4}$ & all  & 4/20 & ? \\
$k_{305}$ & $\prod\cos\phi'_i$ & ? & all & all & 4/20 & \cite{akopyan19_private_pedal_angles} \\
$k_{306}$ & $C_0'$ & ? & all & all & 4/20 & \cite{bialy2020-invariants} \\
$k_{307}$ & $C_2'$ & ? & even & all & 4/20 & ? \\
\hline
\end{tabular}
\caption{Invariants of pedal polygon with respect to the sides of the outer polygon. $^\dagger$ $k_{301}$ means the locus of the outer pedal with respect to a focus is a circle.}
\label{tab:inv-ped-tang}
\end{table}

\subsection{Antipedal Polygons}
\label{sec:antipedal}
\begin{figure}
    \centering
    \includegraphics[width=\textwidth]{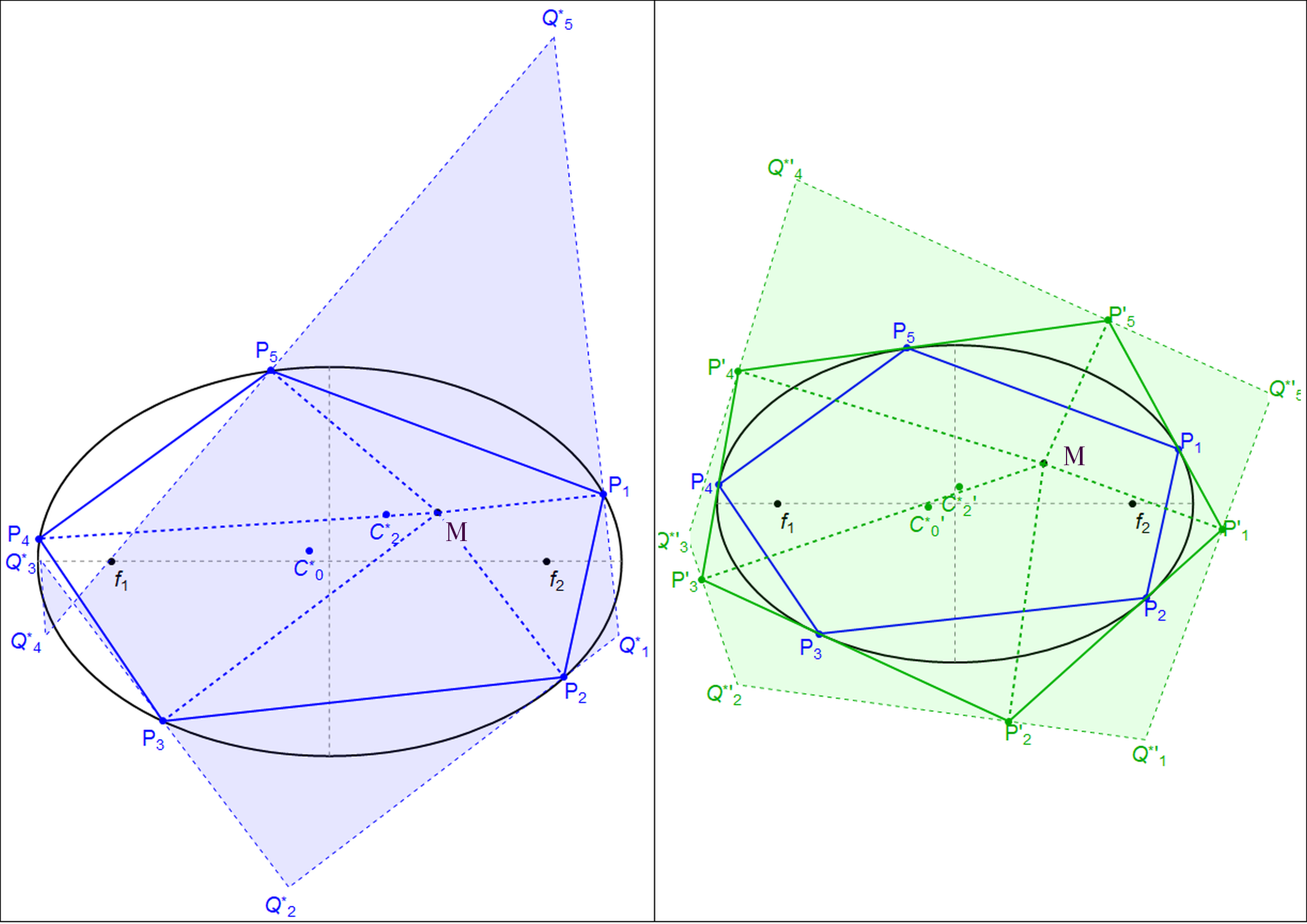}
    \caption{Left (resp. right): Antipedal polygons for $N=5$ from a point $m$ with respect to the $N$-periodic (resp. its outer polygon). Vertex and area centroids $C_0^*,C_2^*$ are also shown.}
    \label{fig:antipedals}
\end{figure}

The antipedal polygons to the $N$-periodic and the outer polygon are shown in Figure~\ref{fig:antipedals}. The antipedal polygon $Q^*_i$ of $P_i$ with respect to $M$ is defined by the intersections of rays shot from every $P_i$ along $(P_i-M)^\perp$. 

%Let $p_i=|P_i-M|$, $p_i'=|P_i'-M|$.

Let $A_m$ denote the area of the $Q^*_i$ polygon and $C_0^*,C_2^*$ its vertex- and signed \footnote{Antipedals can be self-intersecting.} area-centroids. ${C'_0}^*,{C_2'}^*$ refer to centers of antipedals of the outer polygon. Table~\ref{tab:inv-antiped-tang} lists invariants found so far for these polygons.

\begin{table}
\begin{tabular}{|c|l|c|l|c|l|c|}
\hline
code & invariant & value & which N & M & date & proven \\
\hline
$k_{401}$ & $A'\,A_m^*$ & ? & ${\equiv}\,2 \Mod{4}$ & all & 4/20 & ? \\
$k_{402}$ & $A'/A_m^*$ & ? & ${\equiv}\,0 \Mod{4}$ & all & 4/20 & ? \\

$k_{403,a}$ & $A_m\,A_m^*$ & ? & odd & O & 4/20 & ? \\
$k_{403,b}$ & $A_m\,A_m^*$ & ? & ${\equiv}\,0 \Mod{4}$ & $f_1,f_2$ & 4/20 & ? \\
$k_{404}$ & $A_m^*/A_m$ & ? & ${\equiv}\,2 \Mod{4}$ & $f_1,f_2$ & 4/20 & ? \\
$k_{405}$ & $C_0^*$ & ? & even &  O,$f_1,f_2$ & 4/20 & ? \\
$k_{406,a}$ & ${C_0^*}',{C_2^*}'$ & O & even &  O & 4/20 & ? \\
$k_{406,b}$ & ${C_0^*}',{C_2^*}'$ & ? & 4 & $f_1,f_2$ & 4/20 & ? \\
$k_{407}$ & ${C_0^*}'$ & ? & even &  $f_1,f_2$ & 4/20 & ? \\
\hline
\end{tabular}
\caption{Invariants of antipedal polygons.}
\label{tab:inv-antiped-tang}
\end{table}

\subsection{Pedals of Steiner Curvature Centroids}
\label{sec:steiner}
\begin{figure}
    \centering
    \includegraphics{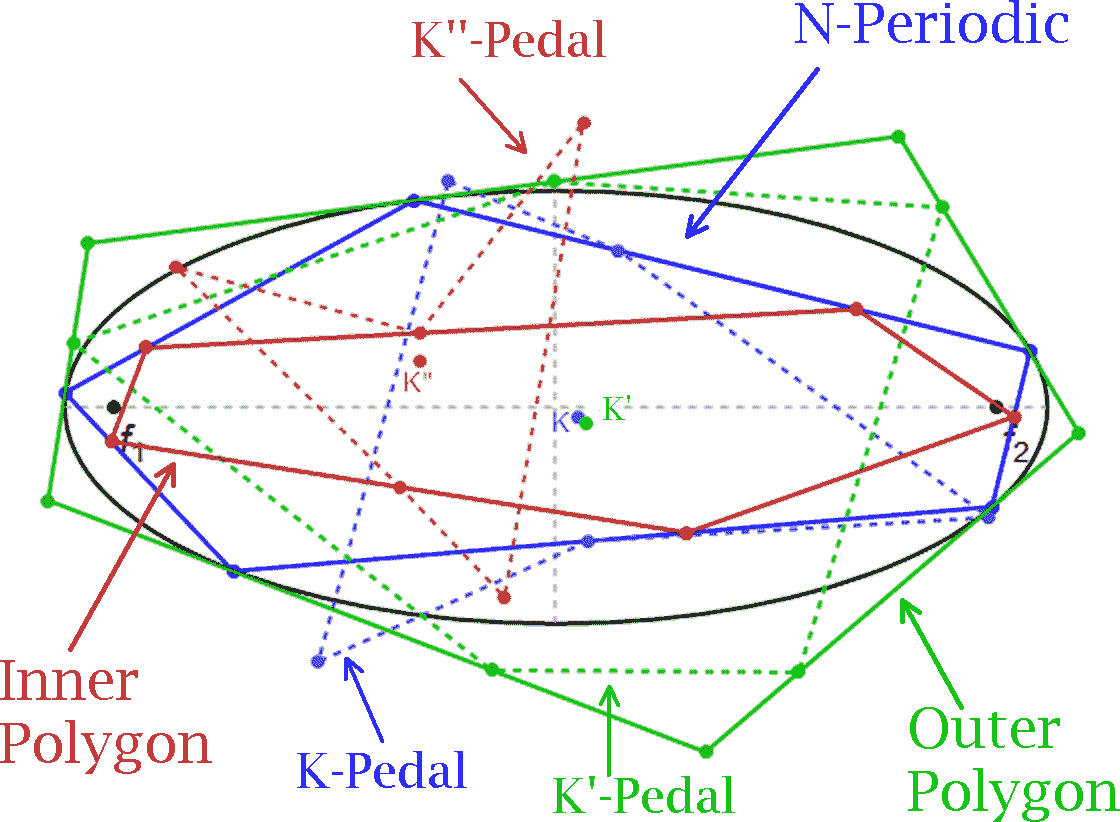}
    \caption{An N-periodic $P$ is shown along with its outer $P'$ and inner $P''$ polygons. Also shown are their Steiner centroids of curvature $K,K',K''$ and the the pedal polygons $P_k,P_k',P_k''$ with respect to said centroids.}
    \label{fig:steiner-centroids}
\end{figure}

Referring to Figure~\ref{fig:steiner-centroids},
let $P,P',P''$ denote as before the N-periodic, outer, and inner polygons, $A, A', A''$ their areas, and $K, K', K''$ their Steiner centroids of curvature \eqref{eqn:steiner}. Let $P_k,P_k',P_k''$ denote the pedal polygons of $P,P',P''$ with respect to $K, K',K''$, and $A_k,A_k',A_k''$ their areas.

When $N$ even, the curvature centroids are stationary at the origin, so invariants described before involving $A,A_m$ (and primed quantities) for $M=O$ apply. For odd $N$, the Curvature Centroids move along individual ellipses concentric with the EB. Invariants are observed appear on Table~\ref{tab:steiner}.

\begin{table}
\begin{tabular}{|c|l|c|c|c|l|}
\hline
code & invariant & value & which N & date & proven \\
\hline
$k_{501}$ & $A/A_k$ & ? & odd &  7/20 & ? \\
%$k_{501,b}$ & $A/A_k$ & ? & ${\equiv}\,2 \Mod{4}$&  7/20 & ? \\
$k_{502}$ & $A'/A_k'$ & ? & odd &  7/20 & ? \\
%$k_{502,b}$ & $A'/A_k'$ & ? &  ${\equiv}\,0 \Mod{4}$ &  7/20 & ? \\
$k_{503}$ & $A''/A_k''$ & ? & odd &  7/20 & ? \\
\hline
\end{tabular}
\caption{Invariants of pedal polygons of $N$-periodic, outer, and inner polygons, with respect to their Steiner Curvature Centroids.}
\label{tab:steiner}
\end{table}

Combining the above with $k_{103}$ and $k_{106}$ one obtains as corollaries the fact that $A_k/A_k'$, $A_k/A_k''$, and $A_k'/A_k''$ are invariant for odd $N$.

\subsection{Pairs of Focal Pedals and Antipedals}
\label{sec:pairs}
Let $\bar{Q}_{1,i}$ and $\bar{Q}_{2,i}$ be the vertices of the pedal polygon with respect to $f_1$ and $f_2$. Define $q_{1,i}=|\bar{Q}_{1,i}-f_1|$ and $q_{2,i}=|\bar{Q}_{2,i}-f_2|$. Likewise, let $\bar{Q}^*_{1,i}$ and $\bar{Q}^*_{2,i}$ be the vertices of the antipedal polygon with respect to $f_1$ and $f_2$. Define $q_{1,i}^*=|\bar{Q}^*_{1,i}-f_1|$ and $q_{2,i}^*=|\bar{Q}^*_{2,i}-f_2|$.

Let $\bar{A}_1$ (resp. $\bar{A}_2$) denote the area of pedal polygon to $N$-periodics wrt $f_1$ (resp. $f_2$) onto the $N$-periodic, and similarly $\bar{A}_1',\bar{A}_2'$ for the outer polygon focus-pedal. Table~\ref{tab:pairs} list invariants so far detected involving pairs of these quantities.

\begin{table}
\begin{tabular}{|c|l|c|l|c|l|c|}
\hline
code & invariant & value & which N & date & proven \\
\hline
$k_{601}$ & $\sum{q_{1,i}}\sum{q_{2,i}}$ & ? & odd & 4/20 & ? \\
$k_{602}$ & $\prod{q_{1,i}}\prod{q_{2,i}}$ & ? & all & 4/20 & ? \\
$k_{603}$ & $\sum{q_{1,i}^*}/\sum{q_{2,i}^*}$ & 1 & all & 5/20 & ? \\
$k_{604,a}$ & $\bar{A}_1.\bar{A}_2$ & ? & odd & 4/20 & ? \\
$k_{604,b}$ & $\bar{A}_1/\bar{A}_2$ & 1 & even & 4/20 & {\scriptsize symmetry} \\
$k_{605,a}$ & $\bar{A}_1'.\bar{A}_2'$ & ? & odd & 4/20 & ? \\
$k_{605,b}$ & $\bar{A}_1'/\bar{A}_2'$ & 1 & even & 4/20 & {\scriptsize symmetry} \\
$k_{606}$ & $\bar{A}_1/\bar{A}_2=\bar{A}_1'/\bar{A}_2'$ & ? & all & 4/20 & ? \\
$k_{607}$ & $\bar{A}_1^*/\bar{A}_2^*$ & 1 & ${\equiv}\,0 \Mod{4}$ & 10/20 & ? \\
$k_{608}$ & $\bar{A}_1^{'*}/\bar{A}_2^{'*}$ & 1 & even & 10/20 & ? \\
$k_{609}$ & $\bar{A}_1''/\bar{A}_2''$ & 1 & even & 10/20 & ? \\
$k_{610}$ & $\bar{A''}_1^{*}/\bar{A''}_2^{*}$ & 1 & even & 10/20 & ? \\
\hline
\end{tabular}
\caption{Invariants between pairs of pedal polygons defined with respect to the foci.}
\label{tab:pairs}
\end{table}

Note $k_{604,a},k_{604,b}$ can be proven via a symmetry argument, namely, area pair are equal since opposite vertices of an even $N$-periodic are reflections about the origin, as will be the pedal polygons from either focus.

\subsection{Evolute Polygons}
\label{sec:evolute-polys}
\begin{figure}
    \centering
    \includegraphics[width=\textwidth]{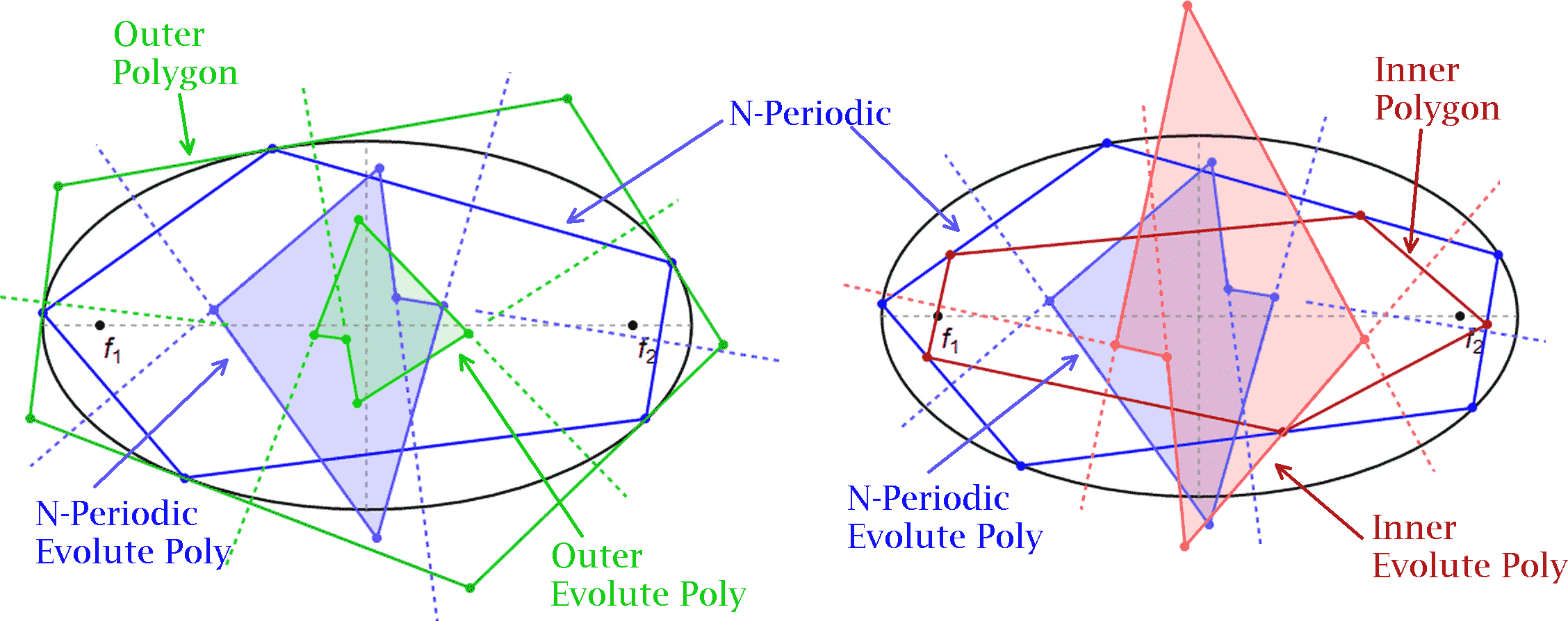}
    \caption{\textbf{Left:} An N-Periodic and its outer polygon are shown along their evolute polygons whose vertices are ordered intersections of perpendicular bisectors. \textbf{Right:} N-periodic,  inner polygon, and their evolute polygons.}
    \label{fig:evolute-polys}
\end{figure}

After \cite{fuchs2017}, let the evolute\footnote{The evolute of a smooth curve is the envelope of the normals \cite[Evolute]{mw}. The perpendicular bisector is its discrete version.} polygon $R_{ev}$ of a generic polygon $R$ have vertices at the intersections of successive pairs of perpendicular bisectors to the sides of $R$; see  Figure~\ref{fig:evolute-polys}. So $P_{ev},P_{ev}',P_{ev}''$ denote the evolute polygons of $P$, $P'$, and $P''$, respectiely, and $A_{ev},A_{ev}',A_{ev}''$ their areas. Trivially, at $N=3$ the latter vanish since perpendicular bisectors concur. At $N=4$, $P'$ is a rectangle, so $A_{ev}'=0$. Area invariants observed for $N>4$ appear on Table~\ref{tab:evolutes}.

\begin{table}
\begin{tabular}{|c|l|c|c|c|c|}
\hline
code & invariant & value & which N & date & proven \\
\hline
$k_{701}$ & $A/A_{ev}$ & ? & $>4$ &  7/20 & ? \\
$k_{702}$ & $A'/A_{ev}'$ & ? & $>4$ &  7/20 & ? \\
$k_{703}$ & $A''/A_{ev}''$ & ? & $>4$ &  7/20 & ? \\
\hline
\end{tabular}
\caption{Area-ratio invariants displayed by the evolute polygons of $N$-periodic, outer, and inner polygons.}
\label{tab:evolutes}
\end{table}

Combining the above with $k_{103}$ and $k_{106}$ one obtains as corollaries the fact that $A_{ev}/A_{ev}'$, $A_{ev}/A_{ev}''$, and $A_{ev}'/A_{ev}''$ are invariant for all $N>4$.

\subsection{Inversive Objects}
\label{sec:inversive}
\begin{figure}
    \centering
    \includegraphics[width=.8\textwidth]{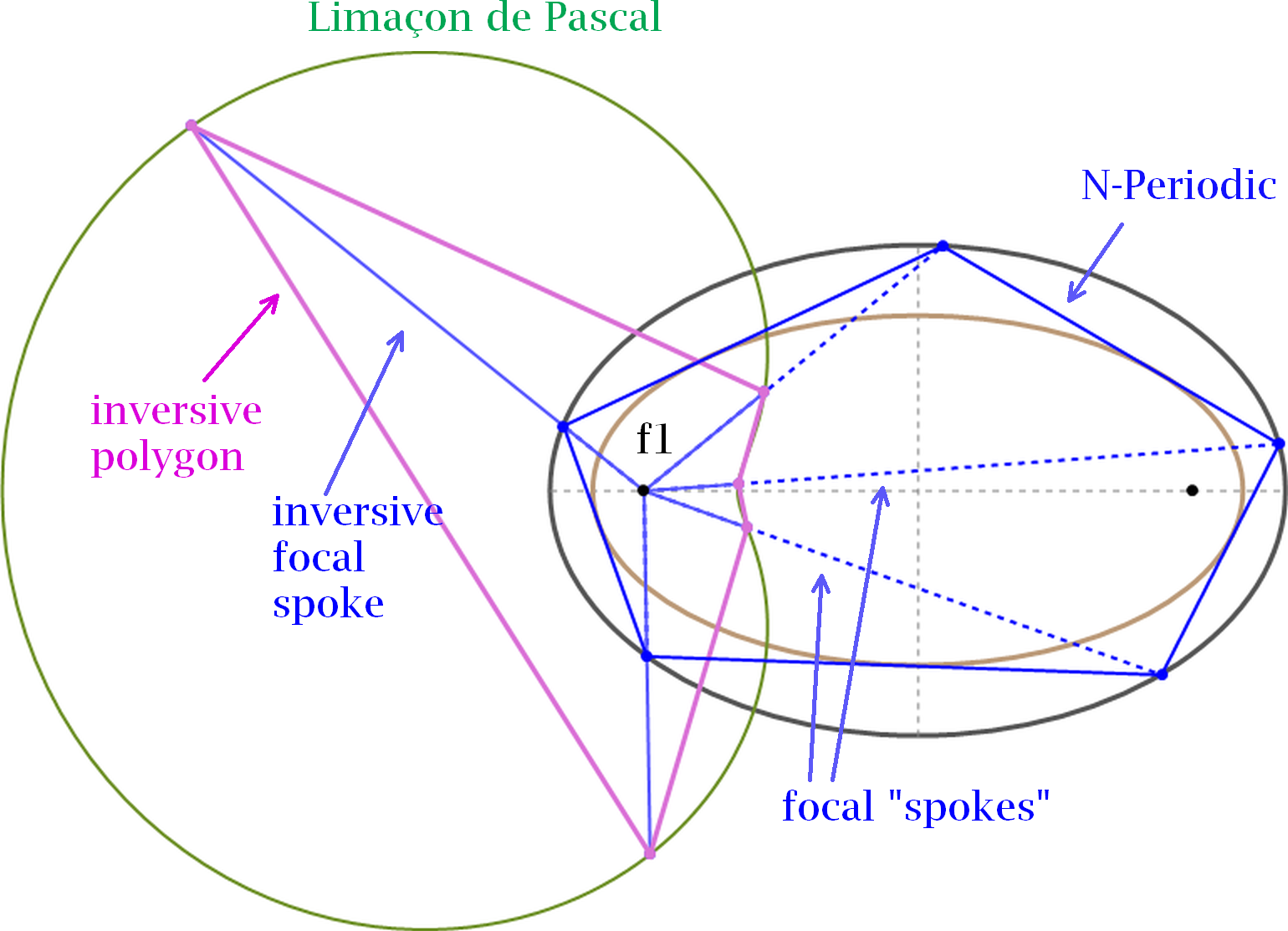}
    \caption{The vertices $P_{j,i}^{-1}$ of the inversive polygon $\P_1$ are obtained by inverting N-periodic vertices with respect to a unit-radius circle (not shown) centered on the left focus. The inversive focal spokes connect said focus to the vertices of $\P_1$.}
    \label{fig:inversive}
\end{figure}

\begin{figure}
    \centering
    \includegraphics[width=.8\textwidth]{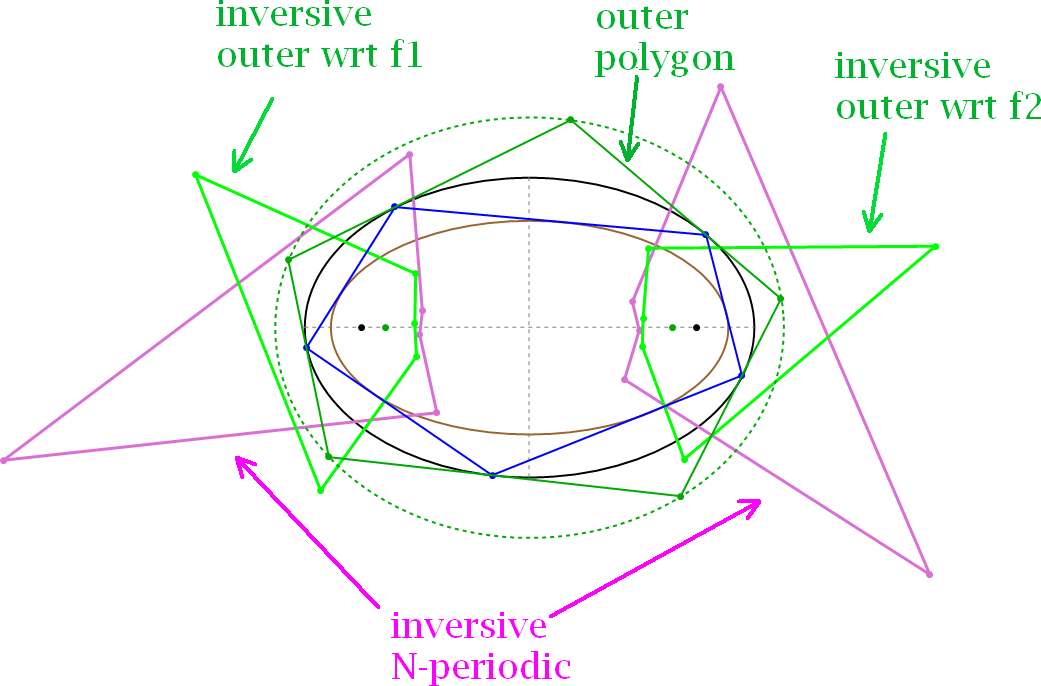}
    \caption{The pair of inversive orbit (resp. outer) polygons is obtaining by inverting the N-periodic (resp. outer polygon) with respect to a circle centered on each focus.}
    \label{fig:inversive-outer}
\end{figure}

Referring to Figure~\ref{fig:inversive}, let $P_{j,i}^{-1}$ denote the inversion of $P_i,i=1,...,N$ with respect to a unit-radius circle centered on focus $f_j$, $j=1,2$, and $d_{j,i}=|P_i-f_j|$. Let $\P_j^\dagger$ denote the polygon with vertices at $P_{j,i}^{-1}$. Let $L_j^\dagger$ denote its perimeter, $A_j^\dagger$ its area, and  $\theta_{j,i}^\dagger$ its ith internal angle. Identical but primed symbols refer to the inversion of the outer polygon (vertices $P_i'$) with respect to $f_j$; see Figure~\ref{fig:inversive-outer}. 

Let $P_i'^\ominus$ (resp. $P_i^\otimes$) denote the inversion of outer (resp. N-periodic) vertices with respect to the billiard (resp. caustic) ellipse. Recall the inversion of a point wrt to an ellipse is the midpoint of the chord joining the tangents from said point to said ellipse \cite{stachel2019-conics}. So the polygon $\P'^\otimes$ (resp. $\P^\ominus$) defined by the $P_i'^\ominus$ (resp. $P_i^\otimes$) has vertices at the side midpoints of N-periodic (resp. inner polygon). Let their areas be denoted $A'^\ominus$ and $A^\otimes$, respectively.

\begin{figure}
    \centering
    \includegraphics[width=\textwidth]{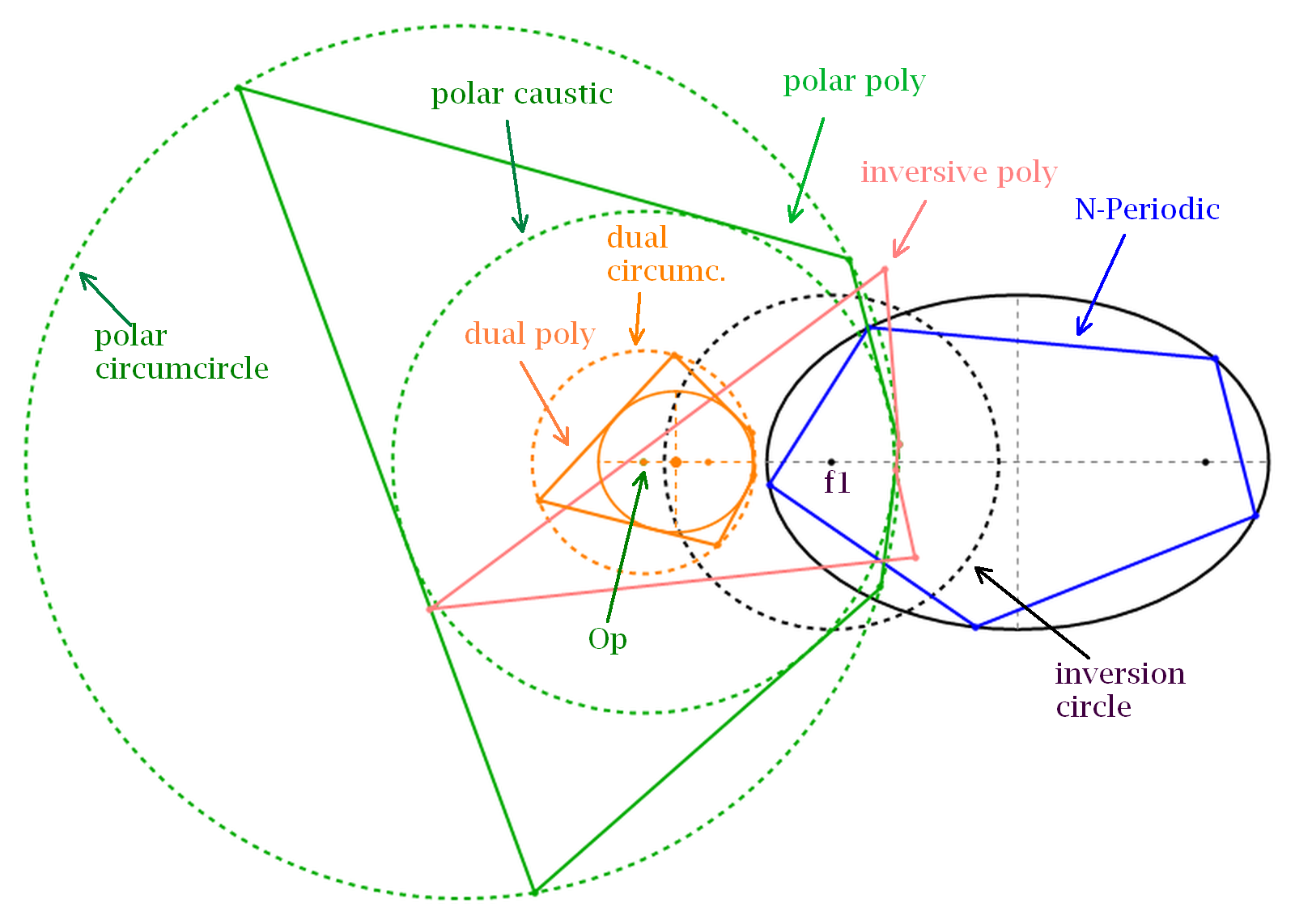}
    \caption{Polygons derived from N-periodics: (i) inversive: vertices are inversions of the $P_i$ wrt to a unit circle centered on a focus, e.g., $f_1$; (ii) polar: antipedal of the inversive wrt to $f_1$. The locus of its vertices is a circle non-concentric with an also circular caustic. Let $O_p$ denote the latter's center. (iii) dual: inversion of the polar with respect to $O_p$. This produces a Poncelet family inscribed in a circle centered on $O_p$ and circumscribed about an ellipse with one focus on $O_p$.}
    \label{fig:polar-dual}
\end{figure}

Referring to Figure~\ref{fig:polar-dual}, let $A_{j,pol}$ and $A_{j,dual}$ denote the areas of the polar and dual polygons with respect to $f_j$, respectively. Let $\psi_{j,i}$ (resp. $w_i$) refers to the polar's ith angle (resp. dual's ith sidelength). Let $A_{j,ant}$ denote the area of the antipedal polygon wrt $f_j$, i.e., $A_m^*$ for $M=f_j$. Recall that the dual is the inverse of the pedal \cite{akopyan2007-conics}, therefore the inversive and antipedal polygons are also inverses of each other. Let $f'_1,f'_2$ be the foci of the elliptic locus of the outer polygons' vertices (non-confocal though concentric and axis aligned with the billiard pair). Let $P_j^{'\ddagger}$ denote the inversion of the $P'_i$ wrt $f'_j$ and $A_j^{'\ddagger}$ denote the of the polygon defined by the $P_j^{'\ddagger}$.

Table~\ref{tab:inversive} lists invariants of inversive objects defined with respect to a chosen focus (e.g. $f_1$), whereas Table~\ref{tab:inversive-pairs} presents invariants involving a pair of inversive objects, defined with respect of both foci. 

\begin{table}
\begin{tabular}{|c|l|c|c|c|c|}
\hline
code & invariant & value & which N & date & proven \\
\hline
$k_{801}$ & $\sum{1/d_{j,i}}$ &  ? & all & 10/20 & \cite{roitman2020-private,akopyan2020-invariants} \\
$k_{802}$ & $L_j^\dagger$ &  ? & all & 10/20 & ? \\
$^{\star\star}{k_{803}}$ & $\sum\cos{\theta_{j,i}^\dagger}$ &  ? & ${\neq}4$ & 10/20 & ? \\
$k_{804,a}$ & $A\,A_j^\dagger$ & ? & ${\equiv}\,0 \Mod{4}$ & 10/20 & ? \\
$k_{804,b}$ & $A\,A_j^\dagger$ & 4 & 4 & 10/20 & ? \\
$k_{805}$ & $A/A_j^\dagger$ & ? & ${\equiv}\,2 \Mod{4}$ & 10/20 & ? \\
$k_{806,a}$ & $A_j'^\dagger/A_j^\dagger$ & ? & all & 10/20 & ?\\
$k_{806,b}$ & $A_j'^\dagger/A_j^\dagger$ & 2 & 4 & 10/20 & ? \\
${^\dagger}k_{807}$ & $A.A^\otimes$ & ? & even & 10/20 & ? \\
$k_{808}$ & $A/A^\otimes$ & ? & odd & 10/20 & ? \\
${^\dagger}k_{809}$ & $A'.{A'^\ominus}$ & ? & even & 10/20 & ? \\
$k_{810}$ & $A'/A'^\ominus$ & ? & odd & 10/20 & ? \\
$^{\ddagger}{k_{811}}$ & $\sum{w^2_{i}}$ &  ? & all & 10/20 & ? \\
${k_{812,a}}$ &
$\sum\cos{\psi_{1,i}}$ & ? & all & 10/20 & ? \\
${k_{812,b}}$ &
$\sum\cos{\psi_{1,i}}$ & 0 & 4 & 10/20 & ? \\
${k_{813}}$ &
$A_{j,pol}/A_{j,inv}$ & ? & all & 10/20 & ? \\
${k_{814}}$ &
$A_{j,pol}/A_{j,dual}$ & ? & all & 10/20 & ? \\
$k_{815}$ & $A_{j,ped}^\dagger.A_{j,dual}$ & ? & odd & 10/20 & ? \\
$k_{816}$ & $A_{j,ped}^\dagger/A_{j,dual}$ & ? & even & 10/20 & ? \\
$k_{817}$ & $A_j^\dagger.A_{j,ant}$ & ? & ${\equiv}\,0 \Mod{4}$ & 10/20 & ? \\
$k_{818}$ & $A_j^\dagger/A_{j,ant}$ & ? & ${\equiv}\,2 \Mod{4}$ & 10/20 & ? \\
\hline
\end{tabular}
\caption{Invariants of inversive objects over the N-periodic family. $^{\star\star}k_{803}$ Co-discovered with Pedro Roitman \cite{roitman2020-private}.  ${^\dagger}k_{807-810}$ The inverted vertices $P_i^\otimes$ (resp. $P_i'^\ominus$) are at the midpoints of the inner (resp. N-periodic orbit) segments. ${^\ddagger}k_{814}$ Discovered by A. Akopyan \cite{akopyan2020-private-polar-sidelengths}.}
\label{tab:inversive}
\end{table}

\begin{table}
\begin{tabular}{|c|l|c|c|c|c|}
\hline
code & invariant & value & which N & date & proven \\
\hline
$k_{901}$ & $\sum{d_{1,i}^{-1}}/\sum{d_{2,i}^{-1}}$ & 1 & all & 10/20 & {\scriptsize from $k_{801}$} \\
$^\star{k_{902}}$ & $\sum{1/(d_{1,i}\,d_{2,i})}$ &  $L/[2 J (a b)^2]$ & all & 10/20 & \cite{stachel2020-private-curvature} \\
$k_{903,a}$ & $A_1^\dagger.A_2^\dagger$ & ? & odd & 10/20 & ? \\
$k_{903,b}$ & $A_1^\dagger/A_2^\dagger$ & 1 & even & 10/20 & {\scriptsize symmetry} \\
$k_{904,a}$ & $A_1'^\dagger.A_2'^\dagger$ & ? & odd & 10/20 & ? \\
$k_{904,b}$ & $A_1'^\dagger/A_2'^\dagger$ & 1 & even & 10/20 & {\scriptsize symmetry} \\
$k_{905}$ & $A_1''^\dagger/A_2''^\dagger$ & 1 & even & 10/20 & ? \\
$k_{906}$ & $A_1'^\ddagger/A_2'^\ddagger$ & 1 & even & 10/20 & ? \\
$k_{907,a}$ & $A_{1,dual}^\dagger.A_{2,dual}$ & ? & odd & 10/20 & ? \\
$k_{907,b}$ & $A_{1,dual}^\dagger/A_{2,dual}^\dagger$ & 1 & even & 10/20 & ? \\
$k_{908,a}$ & $A_{1,ped}^\dagger/A_{2,ped}^\dagger$ & 1 & even & 10/20 & ? \\
$k_{908,b}$ & $A_{1,ped}^\dagger/A_{2,ped}^\dagger$ & 1 & 3 & 10/20 & ? \\
\hline
\end{tabular}
\caption{Invariants of pairs of inversive objects defined with respect to the foci $f_1,f_2$. As observed by A. Akopyan, ${^\star}k_{902}$ is in fact equivalent to $k_{119}$, see \eqref{eqn:curv}.}
\label{tab:inversive-pairs}
\end{table}

\section{Experimental Method}
\label{sec:method}
\begin{figure}
    \centering
    \includegraphics[width=\textwidth]{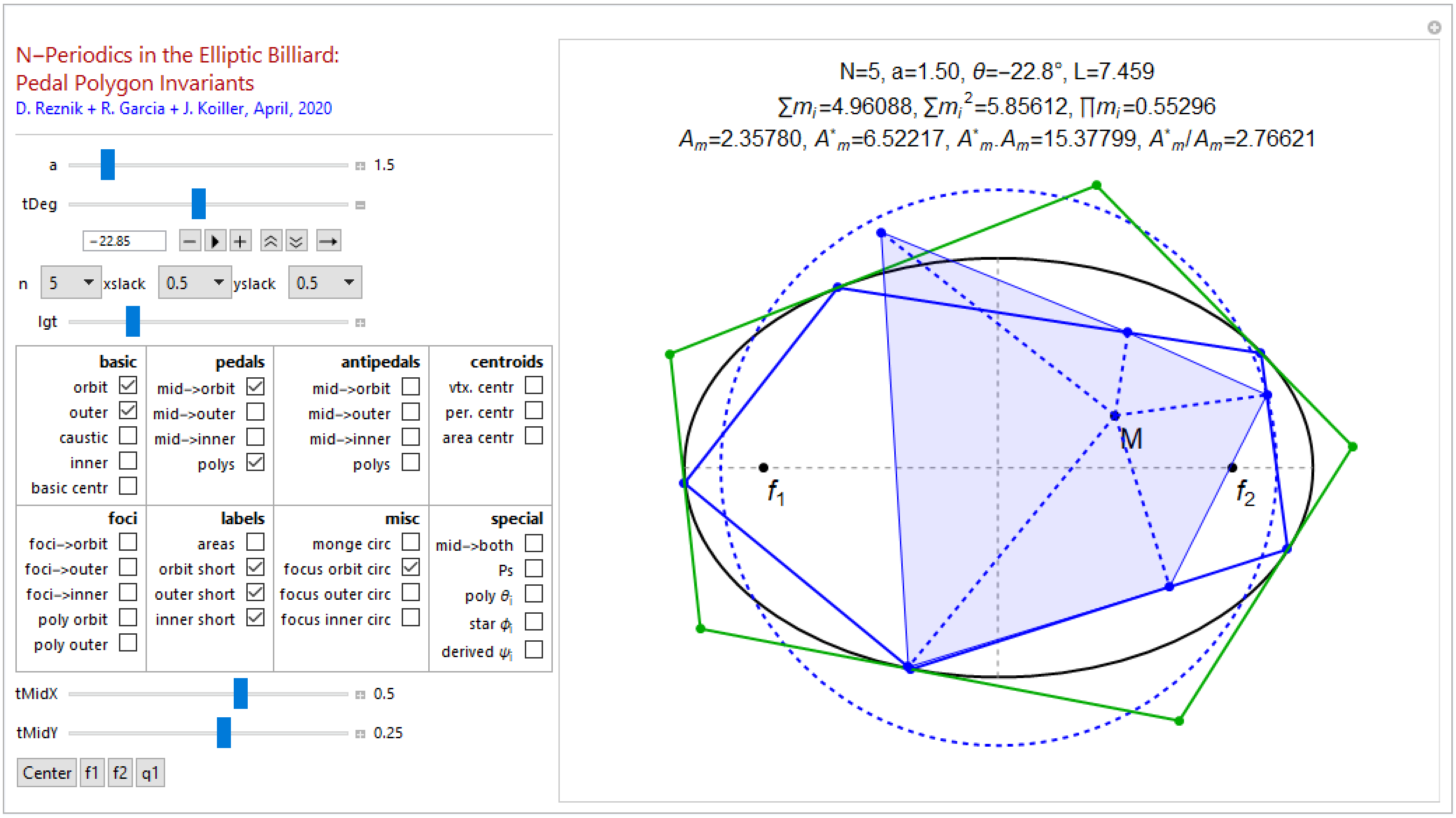}
    \caption{Interactive toolbox written in Wolfram Mathematica \cite{mathematica_v10}. The area on the left permits selection of specific geometries, whereas on the right, the EB, the N-Periodic and derived polygons is displayed. See Videos on Table~\ref{tab:videos}.}
    \label{fig:method}
\end{figure}

A numeric/visualization toolbox was developed in Wolfram Mathematica \cite{mathematica_v10} to accurately calculate and display $N$-periodics while reporting their areas, angles, etc., and those of some derived objects (pedal and inversive polygons, etc.); see Figure~\ref{fig:method}.

Since all trajectories in the Billiard family are tangent to the same confocal caustic, a crucial calculation is to obtain one caustic semiaxis, e.g., $a''$ for a given choice of $a,b$, and $N$. We achieve this for all $N$ via non-linear least-squares minimization \cite{nocedal2006-numopt} of the bisection error. Namely:

\begin{itemize}
    \item Initialize $N$ vertices $P_i$ evenly across the ellipse (pick $t_i, i=1,\ldots,N$ for each), and let $P_1=(a,0)$.
    \item Let $b_i$ be the unit bisector of the $N$-gon sides incident at $P_i$. Let $n_i$ denote the ellipse normal at the $P_i$. The $P_i$ will be a legitimate closed billiard trajectory if all bisectors are perfectly aligned with the local normals, i.e., if $P_i^*$ can be found which make the following error vanish:
    
    \[ \mathcal{E}=\sum_{i=1}^{N}{(n_i^T.b_i)^2} \]
    \item Obtain the unique confocal ellipse tangent to $[a,0]{P_2^*}$.
\end{itemize}

Notice only $N/2$ vertices for $N$ odd (resp. $N/4$ for $N$ even) need to be optimized if one exploits the symmetries of odd (resp. even) vertex positions when $P_1=(a,0)$.

In terms of identifying invariants, we look for quantities which over hundreds of configurations of a given $N$-family are statistically constant, maintained over a range of Billiard aspect ratios. 

\section{Video List}
\label{sec:videos}
Videos of some of the above phenomena have been placed on a Youtube \href{https://bit.ly/379mk1I}{playlist} \cite{dsr_playlist_2020_polygons} and are listed individually on  Table~\ref{tab:videos}.

\begin{table}[H]
\small
\begin{tabular}{|c|l|l|l|}
\hline
Id & Title & N & \texttt{youtu.be/...} \\
\hline
01 & \makecell[lt]{Area Invariants of Pedal\\and Antipedal Polygons} &
3 & \href{https://youtu.be/LN623VjeeFQ}{\texttt{LN623VjeeFQ}} \\
02 & \makecell[lt]{Exploring invariants of N-Periodics\\and pedal polygons} & 3--12 & \href{https://youtu.be/2yXbOV7qf7k}{\texttt{2yXbOV7qf7k}} \\	
03 & Centroid Stationarity of Pedal Polygons &
even & \href{https://youtu.be/j_GD_g8aIbg}{\texttt{j\_GD\_g8aIbg}} \\
04 & \makecell[lt]{Equal sum of distances from foci to\\ vertices of Antipedal Polygon} & 
3--6 & \href{https://youtu.be/6F7Y3UKJzdk}{\texttt{6F7Y3UKJzdk}} \\
05 & \makecell[lt]{Concyclic feet of focal pedals and\\product of sums of lengths for odd N} &
5,6 & \href{https://youtu.be/OT-xAdbOp8o}{\texttt{OT-xAdbOp8o}} \\
06 & \makecell[lt]{Invariant altitudes of N-Periodics\\ and outer polygons I} &
3,4 & \href{https://youtu.be/MvZhWbI6iB8}{\texttt{MvZhWbI6iB8}} \\	
07 & \makecell[lt]{Invariant altitudes of N-Periodics\\and outer polygon II} & 5,6 & \href{https://youtu.be/ZMHLmWXeKrM}{\texttt{ZMHLmWXeKrM}} \\
08 & \makecell[lt]{Sum of focal squared altitudes\\to outer polygon} &
3--8 & \href{https://youtu.be/VUtBRzmbOYU}{\texttt{VUtBRzmbOYU}} \\
09 & \makecell[lt]{Sum of square altitudes from\\arbitrary point to outer polygon} &
5 & \href{https://youtu.be/RNmHROZNGj8}{\texttt{RNmHROZNGj8}} \\	
10 & {Area products of focal pedal polygons}
& 5 & \href{https://youtu.be/sw8pJFMV00w}{\texttt{sw8pJFMV00w}} \\	
11 & \makecell[lt]{Area ratios of Pedal Polygons\\to N-Periodic and outer Polygon} & 
5,6 & \href{https://youtu.be/6F7Y3UKJzdk}{\texttt{6F7Y3UKJzdk}} \\
12 & \makecell[lt]{Invariant Area Ratios to Minimum-Area\\Steiner Pedal Polygons} &
5 & \href{https://youtu.be/f0JwRlu7iaY}{\texttt{f0JwRlu7iaY}} \\
13 & \makecell[lt]{N-Periodic Inversive Invariants} &
5 & \href{https://youtu.be/wkstGKq5jOo}{\texttt{wkstGKq5jOo}} \\
14 & \makecell[lt]{N-Periodic Inversive Objects} &
5 & \href{https://youtu.be/bFsehskizls}{bFsehskizls} \\
15 & \makecell[lt]{Odd N-Periodics in the Elliptic Billiard:\\Invariant Area Product of Focus-Inversive Polygons} &
5 & \href{https://youtu.be/bTkbdEPNUOY}{bTkbdEPNUOY} \\
16 & \makecell[lt]{Invariants of Inversive, Polar, and Dual Polygons\\ derived from Billiard N-Periodics} &
5 & \href{https://youtu.be/qyAHOW32NXY}{qyAHOW32NXY} \\
17 & \makecell[lt]{Centers of N-Periodic Inversive Arcs:\\ Bicentric Poncelet Family w/ Invariants} &
5 & \href{https://youtu.be/mXkk_4RYrnU}{mXkk\_4RYrnU} \\
18 & \makecell[lt]{$N=6$, $a/b=2$ antipedal polygon\\has zero signed area} &
5 & \href{https://youtu.be/fOAES-CzjNI}{\texttt{fOAES-CzjNI}} \\
\hline
\end{tabular}
\caption{Youtube list of videos about invariants of N-Periodics, the last column provides the link.}
\label{tab:videos}
\end{table}

\section*{Acknowledgments}
We would like to thank Olga Romaskevitch, Sergei Tabachnikov, Richard Schwartz, Arseniy Akopyan, Hellmuth Stachel, Alexey Glutsyuk, Corentin Fierobe, Maxim Arnold, and Pedro Roitman for useful discussions and insights.

The second author is fellow of CNPq and coordinator of Project PRONEX/CNPq/FAPEG 2017 10 26 7000 508.

\appendix
\section{Table of Symbols}
\label{app:symbols}
\begin{table}[H]
\begin{tabular}{|c|l|l|}
\hline
symbol & meaning \\
\hline
$O,f_1,f_2,N$ & center and foci of billiard \\
$(a,b),(a'',b'')$ & billiard (resp. caustic) major, minor semi-axes \\
$N,L,J$ & trajectory sides, inv. perimeter and Joachimsthal's constant \\
$P_i,P_i',P_i''$ &
$N$-periodic, outer, inner polygon vertices \\
$d_{j,i}$, $l_i$, $r_i$ & distances $|P_i-f_j|$, $|P''_i-P_i|$, $|P_{i+1}-P''_i|$ \\
$\theta_i,\theta_i',\alpha_{j,i}$ & $N$-periodic, outer poly, and $P_i f_j P_{i+1}$ angles \\
$A,A',A''$ &  $N$-periodic, outer, inner areas \\
\hline
$M$ & a point in the plane of the billiard \\
$Q_i,Q_i',Q_i''$ &
\makecell[tl]{vertices of the pedal polygon of\\$N$-periodic,outer,inner polygon wrt $M$} \\
$Q_i^*,{Q_i^*}',{Q_i^*}''$ & \makecell[lt]{vertices of the antipedal polygon of the\\$N$-periodic,outer,inner polygon wrt $M$}\\
$\phi_i,\phi'_i,\phi''_i$ & \makecell[lt]{ith angle of pedal polygon of\\$N$-periodic,outer,inner polygons wrt $M$} \\
$A_m,A_m',A_m^*$ & areas of $Q_i,Q_i',Q_i^*$ polygons \\ 
%$p_i,p_i'$ & $|P_i-M|$, $|P_i'-M|$ \\
%$d_{0,i},d'_{0,i}$ & $|Q_i-O|$, $|Q_i'-O|$ \\
\hline
$C_0,C_0',C_0^*$ & vertex centroids of the $Q_i,Q_i',Q_i^*$ polygons \\
$C_2,C_2',C_2^*$ & area centroids of the $Q_i,Q_i',Q_i^*$ polygons \\
${C_0^*}',{C_2^*}'$ & vertex, area centroids of the ${Q_i^*}'$ polygon \\
\hline
$\bar{Q}_{j,i},\bar{Q}^*_{j,i}$ & vertices of N-periodic pedal, antipedal polygon wrt. $f_j$ \\
${q}_{j,i},{q}^*_{j,i}$ & $|\bar{Q}_{j,i}-f_j|$ and  $|\bar{Q}^*_{j,i}-f_j|$\\
$\bar{A}_j,\bar{A}_j^*$ & \makecell[lt]{areas of $\bar{Q}_{j,i}$ and $\bar{Q}_{j,i}^*$ polygons} \\
\hline
$K,K',K''$ & Steiner centroids of curvature of $P,P',P''$ \\
$P_k,P_k',P_k''$ & Pedal Polygons of $P,P',P''$ wrt. $K,K',K''$ \\
$A_k,A_k',A_k''$ & Areas of $P_k,P_k',P_k''$ \\
$P_{ev},P_{ev}',P_{ev}''$ & Evolute Polygons of $P,P',P''$ \\
$A_{ev},A_{ev}',A_{ev}''$ & Areas of $P_{ev},P_{ev}',P_{ev}''$ \\
\hline
$P_{j,i}^{-1}$ & \makecell[lt]{inversion of $P_i$ wrt. to unit-radius circle centered on $f_j$} \\
$\P_j^\dagger,L_j^\dagger,A_j^\dagger$ & polygon defined by $P_{j,i}^{-1}$, its perimeter, and area \\
$\theta_{j,i}^\dagger$ & internal angle of $\P_j^\dagger$ at $P_{j,i}^{-1}$ \\
$f_1',f_2'$ & foci of the elliptic locus of the $P_i'$ \\
$P'_j{^\ddagger},A'_j{^\ddagger}$ & inversion of $P_i'$ wrt to $f_j'$, area of $P'_j{^\ddagger}$ poly\\
$\P^\otimes,\P'^\ominus$ & inversion of $P_i$ (resp. $P_i'$) wrt caustic (resp. billiard) \\
$A^\otimes,A'^\ominus$ & areas of $\P^\otimes,\P'^\ominus$ \\
\hline
$\psi_{j,i},w_i$ & ith angle of polar (sidelength of dual) polygon wrt $f_j$\\
$A_{j,pol},A_{j,dual}$ & area of polar, dual polygon wrt $f_j$ \\
$A_{j,ped},A_{j,ant}$ & area of pedal, antipedal polygon wrt $f_j$ \\
\hline
\end{tabular}
\caption{Symbols used in the invariants. Note $i=1,...,N$ and $j=1,2$. Note any single (resp. double) primed quantities apply to the outer (resp. inner) polygons.}
\label{tab:symbols}
\end{table}

\bibliographystyle{maa}
\bibliography{references} 

\begin{thebibliography}{10}
\expandafter\ifx\csname urlstyle\endcsname\relax
 \providecommand{\url}[1]{doi:\discretionary{}{}{}#1}\else
 \providecommand{\url}{doi:\discretionary{}{}{}\begingroup
  \urlstyle{rm}\Url}\fi

\bibitem{akopyan19_private_pedal_angles}
Akopyan, A. (2020).
\newblock Angles $\phi=\pi-\theta_i$ (resp. $\phi'=\phi-\theta_i'$), so
  equivalent to invariant sum (resp. product) of cosines.
\newblock Private Communication.

\bibitem{akopyan19_private_half_sines}
Akopyan, A. (2020).
\newblock {Corollary of Theorem 6 in {A}kopyan et al., ``Billiards in Ellipses
  Revisited'' (2020)}.
\newblock Private Communication.

\bibitem{akopyan2020-private-affine}
Akopyan, A. (2020).
\newblock Follows from previous results: the construction is affine and holds
  for any two concentric conics.
\newblock Private Communication.

\bibitem{akopyan20_private_cyclic_feet}
Akopyan, A. (2020).
\newblock Perpendicular feet to {N}-periodic or its tangential polygon are
  cyclic.
\newblock Private Communication.

\bibitem{akopyan2020-private-polar-sidelengths}
Akopyan, A. (2020).
\newblock Sum of squared sidelengths of focus-polar polygon is invariant.
\newblock Private Communication.

\bibitem{akopyan2020-invariants}
Akopyan, A., Schwartz, R., Tabachnikov, S. (2020).
\newblock Billiards in ellipses revisited.
\newblock \emph{Eur. J. Math.}
\newblock \url{doi.org/10.1007/s40879-020-00426-9}.

\bibitem{akopyan2007-conics}
Akopyan, A., Zaslavski, A. (2007).
\newblock \emph{Geometry of Conics}.
\newblock Mathematical world. Providence: American Mathematical Society.

\bibitem{fuchs2017}
Arnold, M., Fuchs, D., Izmestiev, I., Tabachnikov, S., Tsukerman, E. (2017).
\newblock Iterating evolutes and involutes.
\newblock \emph{Discrete Comput. Geom.}, 58: 80--–143.
\newblock \url{doi.org/10.1007/s00454-017-9890-y}.

\bibitem{bialy2020-invariants}
Bialy, M., Tabachnikov, S. (2020).
\newblock {Dan Reznik's} identities and more.
\newblock \emph{Eur. J. Math.}
\newblock \url{doi.org/10.1007/s40879-020-00428-7}.

\bibitem{chang1988}
Chang, S.-J., Friedberg, R. (1988).
\newblock Elliptical billiards and {P}oncelet’s theorem.
\newblock \emph{J. Math. Phys.}, 29: 1537--1550.

\bibitem{caliz2020-area-product}
Chavez-Caliz, A. (2020).
\newblock More about areas and centers of {Poncelet} polygons.
\newblock \emph{Arnold Math J.}
\newblock \url{doi.org/10.1007/s40598-020-00154-8}.

\bibitem{dragovic11}
Dragovi\'{c}, V., Radnovi\'{c}, M. (2011).
\newblock \emph{Poncelet Porisms and Beyond: Integrable Billiards,
  Hyperelliptic Jacobians and Pencils of Quadrics}.
\newblock Frontiers in Mathematics. Basel: Springer.

\bibitem{dragovic2014-poncelet}
Dragovi\'{c}, V., Radnovi\'{c}, M. (2014).
\newblock Bicentennial of the great {P}oncelet theorem (1813–2013): current
  advances.
\newblock \emph{Bulletin Amer. Math. Soc.}, 51(3): 373--445.

\bibitem{dragovic2006-survey}
Dragović, V., Radnović, M. (2006).
\newblock A survey of the analytical description of periodic elliptical
  billiard trajectories.
\newblock \emph{J. of Math. Sci.}, 135: 3244--3255.
\newblock \url{doi.org/10.1007/s10958-006-0154-2}.

\bibitem{garcia2020-invariants}
Garcia, R., Reznik, D. (2020).
\newblock Explicit expressions for some elliptic billiard invariants.
\newblock In preparation.

\bibitem{garcia2020-new-properties}
Garcia, R., Reznik, D., Koiller, J. (2020).
\newblock New properties of triangular orbits in elliptic billiards.
\newblock \emph{Amer. Math. Monthly}, to appear.

\bibitem{stachel2019-conics}
Glaeser, G., Stachel, H., Odehnal, B. (2016).
\newblock \emph{The Universe of Conics: From the ancient Greeks to 21st century
  developments}.
\newblock Berlin: Springer.

\bibitem{izmestiev2017-ivory}
Izmestiev, I., Tabachnikov, S. (2017).
\newblock Ivory’s theorem revisited.
\newblock \emph{Journal of Integrable Systems}, 2.
\newblock \url{doi.org/10.1093/integr/xyx006}.

\bibitem{kaloshin2018}
Kaloshin, V., Sorrentino, A. (2018).
\newblock On the integrability of {B}irkhoff billiards.
\newblock \emph{Phil. Trans. R. Soc.}, A(376).

\bibitem{nocedal2006-numopt}
Nocedal, J., Wright, S. (2006).
\newblock \emph{Numerical Optimization}, chap.~10.
\newblock New York: Springer, 2nd ed., pp. 245--269.
\newblock Least-Squares Problems.

\bibitem{nomizu1994-affine}
Nomizu, K., Sasaki, T. (1994).
\newblock \emph{Affine Differential Geometry}.
\newblock Cambridge: Cambridge University Press.

\bibitem{preparata1988}
Preparata, F., Shamos, M. (1988).
\newblock \emph{Computational Geometry - An Introduction}.
\newblock New York: Springer-Verlag, 2nd ed.

\bibitem{ros2014-ellipsoids}
Ramírez-Ros, R. (2014).
\newblock On {C}ayley conditions for billiards inside ellipsoids.
\newblock \emph{Nonlinearity}, 27(5): 1003--1028.

\bibitem{dsr_playlist_2020_invariants}
Reznik, D. (2020).
\newblock Playlist for {``Invariants of 3- and 4-Periodics in the Elliptic
  Billiard''}.
\newblock YouTube.
\newblock \url{bit.ly/3aNqgqU}.

\bibitem{dsr_playlist_2020_polygons}
Reznik, D. (2020).
\newblock Playlist for {``Invariants of N-Periodics in the Elliptic
  Billiard''}.
\newblock YouTube.
\newblock \url{bit.ly/2xeVGYw}.

\bibitem{reznik2019-intelligencer}
Reznik, D., Garcia, R., Koiller, J. (2020).
\newblock Can the elliptic billiard still surprise us?
\newblock \emph{Math. Intelligencer}, 42: 6--17.
\newblock \url{rdcu.be/b2cg1}.

\bibitem{roitman2020-private}
Roitman, P. (2020).
\newblock Investigation of n-periodic invariants involving the curvature of the
  ellipse.
\newblock Private Communication.

\bibitem{rozikov2018}
Rozikov, U.~A. (2018).
\newblock \emph{An Introduction To Mathematical Billiards}.
\newblock Singapore: World Scientific Publishing Company.

\bibitem{sergei2016-com}
Schwartz, R., Tabachnikov, S. (2016).
\newblock Centers of mass of {P}oncelet polygons, 200 years after.
\newblock \emph{Math. Intelligencer}, 38(2): 29--34.
\newblock \url{doi.org/10.1007/s00283-016-9622-9}.

\bibitem{stachel2020-private-curvature}
Stachel, H. (2020).
\newblock Closed form expression for $k_{119}$.
\newblock Private Communication.

\bibitem{stachel2020-private-joachimsthal}
Stachel, H. (2020).
\newblock Joachmisthal's constant in terms of $a$, $b$ and $a''$.
\newblock Private Communication.

\bibitem{stachel2020-private-k113-k116}
Stachel, H. (2020).
\newblock Proofs for $k_{113}$ and $k_{116}$.
\newblock Private Communication.

\bibitem{steiner1838}
Steiner, J. (1838).
\newblock {Über den Krümmungs-Schwerpunkt ebener Curven}.
\newblock \emph{{Abhandlungen der Königlichen Akademie der Wissenshaften zu
  Berlin}}: 19--91.

\bibitem{sergei91}
Tabachnikov, S. (2005).
\newblock \emph{Geometry and Billiards}, vol.~30 of \emph{Student Mathematical
  Library}.
\newblock Providence, RI: American Mathematical Society.
\newblock \url{bit.ly/2RV04CK}.

\bibitem{mw}
Weisstein, E. (2019).
\newblock Mathworld.
\newblock \url{mathworld.wolfram.com}.

\bibitem{mathematica_v10}
Wolfram, S. (2019).
\newblock Mathematica, version 10.0.

\end{thebibliography}

\end{document}